\documentstyle[12pt]{article}
\renewcommand{\b}{\beta}

\newfont{\bbb}{msbm10 scaled \magstep1}

\newcommand{\link}{\!:\!}

\newcommand{\binom}[2]{{#1}\choose{#2}}
\newcommand{\socle}{\mbox{\rm{socle}}\,}

\def\demo{\noindent{\em Proof. }}
\def\QED{\hfill$\Box$}

\newcommand{\be}{\begin{equation}}
\newcommand{\ee}{\end{equation}}
\newcommand{\ba}{\begin{eqnarray}}
\newcommand{\ea}{\end{eqnarray}}
\newtheorem{theorem}{Theorem}[section]

\newtheorem{remark}[theorem]{Remark}

\newtheorem{example}[theorem]{Example}
\newtheorem{conjecture}[theorem]{Conjecture}

\begin{document}

\title{\large\sc\bf On the Betti numbers of some \\ Gorenstein ideals}
\author{
{\normalsize\sc Matthew Miller\thanks{\rm Research supported by the National Science Foundation}}\vspace{-0.75mm}\\
{\small\em Department of Mathematics,  University of South Carolina }\vspace{-1.75mm}\\
{\small\em Columbia, SC 29208 USA}
  \and
{\normalsize\sc Rafael H. Villarreal\thanks{\rm Partially supported by COFAA--IPN, CONACyT  and SNI--SEP}}\vspace{-0.75mm}\\
{\small\em Depto. de Mat., ESFM del Instituto Polit\'ecnico Nacional}\vspace{-1.75mm}\\
{\small\em Unidad Adolfo L\'opez Mateos}\vspace{-1.75mm}\\
{\small\em 07300 M\'exico, D.F. MEXICO}
}
\maketitle

\begin{abstract} Assume $R$ is a polynomial ring over a field and 
$I$ is a homogeneous Gorenstein ideal of codimension $g\ge3$ and
initial degree $p\ge2$. We prove that the number of minimal generators $\nu(I_p)$ of $I$ that are in degree $p$ is bounded above by
$\nu_0={p+g-1\choose g-1}-{p+g-3\choose g-1}$, which is the number of minimal generators of the defining ideal of the extremal Gorenstein algebra of codimension $g$ and initial degree $p$. Further, $I$ is itself extremal if $\nu(I_p)=\nu_0$. 
\end{abstract}

%
%
\section{Introduction}

\hskip1\parindent Assume $R$ is a polynomial ring over a field and 
$I$ is an homogeneous Gorenstein ideal of codimension $g\ge3$ and
initial degree $p\ge2$. We have the following conjectures on the minimal number of generators of the ideal generated by the forms of degree $p$ in $I$. 

\begin{conjecture}\label{conj:1}\sl
\begin{enumerate}
\item Always $\nu(I_p)\leq \nu_0={p+g-1\choose g-1}-{p+g-3\choose g-1}$, and only certain values 
of $\nu(I_p)$ are possible.
\item If $\nu(I_p)=\nu_0$ then
$I$ is extremal in the sense of \cite{S}, or equivalently,
$$
e(R/I)= e(g,p)={{g+p-1}\choose{g}}+{{g+p-2}\choose g}.
$$
Consequently if $\nu(I_p)=\nu_0$ then $I=(I_p)$. 
\end{enumerate}
\end{conjecture}

These  estimates involve comparisons with the
numerical invariants of Schenzel's \cite{S} extremal Gorenstein
algebras.
If $I$ is a graded Gorenstein ideal of codimension $g$ and initial
degree $p$, then  a
consequence of the Macaulay-Stanley characterization \cite{St} of the
Hilbert function of $R/I$ is that the multiplicity of $R/I$ satisfies
$$
e(R/I)\ge e(g,p)=
{{g+p-1}\choose{g}}+{{g+p-2}\choose g}.
$$
(A version of this estimate for non-graded Gorenstein algebras 
is apparently an open problem; see \cite{R-V}.)
Given codimension $g$ and initial degree $p$, a graded
Gorenstein algebra $R/J$ with  multiplicity $e(R/J)=e(g,p)$ must have a pure  
almost linear minimal resolution (in particular $J=(J_p)$ and  $\nu_0$ is the number of generators of $J$), 
and there are similar formulas for the other betti numbers of $R/J$. 
Conversely, if the resolution of $R/J$ is pure and almost linear, or
equivalently, $R/J$ is extremal, then
$J=(J_p)$, $\nu(J)=\nu_0$, 
and $e(R/J)=e(g,p)$.

Various results in the literature dealing with Cohen-Macaulay ideals
 (such as \cite{E,E-R-V,R-V})
give upper bounds for $\nu(I)$ in terms of codimension, initial
degree, {\em and multiplicity } $e(R/I)$. Part of the intent of our
conjectures is to elucidate the multiplicity information that is
already determined (if it is at all!) by the codimension and initial
degree. In general this sounds quite implausible, but for Gorenstein
algebras it looks promising that something can be said along these
lines. Since, however, symmetry of the $h$-vector $H(R/I)$ appears to play a
major role, one would not expect any such results to generalize to the
non-graded case. 

One might optimistically hope for even stronger estimates than suggested by 
Conjecture \ref{conj:1}, for instance that $\nu(I)\le\nu_0$ or even $\b_i(R/I)\le\b_i(R/J)$ for all $i$ (that is remove the restriction to the degree $p$ generators of $I$, and then pass from the first betti number of $R/I$ to the entire minimal free resolution). We have produced a considerable amount of computational evidence by using the program MACAULAY, but the same program also enabled us to find a counterexample. We are now looking for reasonable side conditions under which these stronger estimates might hold.

\section{Hilbert function techniques}

\begin{theorem}\label{main}If $I$ is a graded Gorenstein ideal of
codimension $g\ge3$ and initial degree $p\ge2$, then
\[ \nu(I_p)\leq \nu_0={p+g-1\choose g-1}-{p+g-3\choose g-1}, \]
and $I$ is itself extremal if equality holds.
\end{theorem}

\demo If either $\nu(I_p)>\nu_0$, or $\nu(I_p)=\nu_0$ and $I$ is
not extremal, then by the symmetry of the $h$-vector $H(R/I)$ there
is some $j\ge p$ so that
$$H(R/I)=
(h_0,h_1,\ldots, h_{p-1},h_p,
\ldots,h_j,h_{p-1} ,\ldots, h_1, h_0),$$
where $h_j=h_p\leq {p+g-3\choose g-1}=h_{p-2}$. The idea of
the argument is to use the Macaulay estimate
(see \cite{E,E-R-V,R,St}) for $h_{j+1}$ in terms of $j$ and $h_j$
to see that such a small value of $h_j$ can not grow to
such a large value of $h_{j+1}=h_{p-1}={p+g-2\choose g-1}$. We recall that this estimate
is calculated from the binomial expansion for $h_j$:
\be
h_j={{a_j}\choose j}+{{a_{j-1}}\choose{j-1}}
+\cdots+{{a_i}\choose i}\label{eq:h_j}
\ee
where $a_j>a_{j-1}>\cdots> a_i\ge i\ge1$. Then
\be
h_{j+1}\leq {{a_j+1}\choose{j+1}}+ {{a_{j-1}+1}\choose{j}}
+\cdots+{{a_i+1} \choose{i+1}}.\label{eq:h_j+1}
\ee
 We may assume that $h_j>j$, for if
not, then $h_j\leq j$ would imply that $a_{\ell}=\ell$ for all $\ell$, and hence  $h_j \geq h_{j+1}$, which  contradicts our assumption. Notice that 
by grouping the terms of (\ref{eq:h_j}) according to the value of $a_{\ell}-\ell$ the binomial
expansion for $h_j$ can be written as
\ba
h_j &= & \sum_{n=0}^{r}\left[ {j_n+k_n\choose j_n}
+{j_n-1+k_n\choose j_n-1}+\cdots+{j_n-i_n+k_n\choose j_n-i_n}\right] \nonumber\\
&= & \sum_{n=0}^{r}\left[{j_n+k_n+1\choose
j_n}-{j_n-i_n+k_n\choose j_n-i_n-1}\right],
\label{eq:h_j=expanded}
\ea
where $a_j-j=k_0>k_1>\cdots>k_r\geq 0$,
$j=j_0>j_1>\cdots>j_r$, $j_r-i_r=i$, and  $j_n=j_{n-1}-i_{n-1}-1$ for
$1\leq n\leq r$. Set $k=k_0$.
Since $ p\le j$ and
${j+k\choose j}\leq h_j\le{p+g-3\choose g-1}$ it follows that $k\leq g-2$.
From (\ref{eq:h_j}), (\ref{eq:h_j+1}), and (\ref{eq:h_j=expanded}), together with Pascal's identity and ${a+b+1\choose b+1}=\frac{a+1}{b+1}{a+b+1\choose b}$,  we have
\begin{eqnarray*}
h_{j+1}-h_j &\leq &\sum_{n=0}^{r}\left[{j_n+k_n+1\choose j_n+1}-
{j_n-i_n+k_n\choose j_n-i_n}\right]\\
& = &\sum_{n=0}^r\left[\frac{k_n+1}{j_n+1}
{j_n+k_n+1\choose j_n}-\frac{k_n+1}{j_n-i_n}{j_n-i_n+k_n\choose
j_n-i_n-1}\right].
\end{eqnarray*}
On the other hand from the upper bound on $h_j$ and $h_{j+1}=h_{p-1}$ we see that
\[
\frac{g-1}{p-1} h_j\leq {p+g-3\choose g-2}
\leq h_{j+1}-h_j.
\]
Since $(g-1)/(p-1)>(k+1)/(j+1)$ it follows from (\ref{eq:h_j=expanded}) and the last two inequalities that 
\ba
F_0 &= &\sum_{n=0}^r\left[\left(\frac{k+1}{j+1}-\frac{k_n+1}{j_n+1}\right)
{j_n+k_n+1\choose j_n}\right] \nonumber \\
&  & \quad \mbox{} +\sum_{n=0}^r\left[\left(\frac{k_n+1}{j_n-i_n}-\frac{k+1}{j+1}\right)
{j_n-i_n+k_n\choose j_n-i_n-1}\right] \quad < \quad 0.\label{eq:F_0<0}
\ea
For $0\leq s\leq r$ we set
\begin{eqnarray*}
F_s &= &\sum_{n=s}^r\left[\left(\frac{k_s+1}{j_s+1}-\frac{k_n+1}{j_n+1}\right)
{j_n+k_n+1\choose j_n}\right]\\
 &  & \quad \mbox{} +\sum_{n=s}^r\left[\left(\frac{k_n+1}{j_n-i_n}-\frac{k_s+1}{j_s+1}\right)
{j_n-i_n+k_n\choose j_n-i_n-1}\right].
\end{eqnarray*}
To derive a contradiction we are going to show the following
inequalities
\[
 F_0>F_1>\cdots >F_r=\left(\frac{k_r+1}{j_r-i_r}-
\frac{k_r+1}{j_r+1}\right){a_i\choose i-1}>0.
\]
Assume $1\leq s+1\leq r$. Notice that
\begin{eqnarray*}
0 &\leq &\sum_{n=s}^{r-1}\left[{j_n-i_n+k_n\choose j_n-i_n-1}-
{j_{n+1}+k_{n+1}+1\choose j_{n+1}}\right]
+{j_r-i_r+k_r\choose j_r-i_r-1}\\
&= &{j_s-i_s+k_s\choose j_s-i_s-1}-\sum_{n=s+1}^r
\left[{j_n+k_n+1\choose j_n}-{j_n-i_n+k_n\choose j_n-i_n-1}\right].
\end{eqnarray*}
Therefore
\begin{eqnarray*}
{j_s-i_s+k_s\choose j_s-i_s-1}\geq\sum_{n=s+1}^r
{j_n+k_n+1\choose j_n}-\sum_{n=s+1}^r{j_n-i_n+k_n\choose j_n-i_n-1}.
\end{eqnarray*}
Let $A_s $ denote the first summation and $B_s $ the second in
 this last inequality; clearly $A_s-B_s>0$. Also note that 
\[
\frac{k_s+1}{j_s-i_s}>\frac{k_s+1}{j_s+1}\qquad\mbox{and}\qquad \frac{k_s+1}{j_s-i_s}>\frac{k_{s+1}+1}{j_{s+1}+1}.
\]
Putting these together, we compute 
\begin{eqnarray*}
F_s-F_{s+1}&= &\left(\frac{k_s+1}{j_s-i_s}-\frac{k_s+1}{j_s+1}\right)
{j_s-i_s+k_s\choose j_s-i_s-1}\\
 &  &  \qquad \mbox{}+\left(\frac{k_s+1}{j_s+1}-\frac{k_{s+1}+1}{j_{s+1}+1}\right)A_s 
+\left(\frac{k_{s+1}+1}{j_{s+1}+1}-\frac{k_s+1}{j_s+1}\right)B_s \\
\pagebreak[2]
&\geq & \left(\frac{k_s+1}{j_s-i_s}-\frac{k_s+1}{j_s+1}\right)(A_s -B_s )+
\left(\frac{k_s+1}{j_s+1}-\frac{k_{s+1}+1}{j_{s+1}+1}\right)A_s \\
 &  &  \qquad \mbox{} +\left(\frac{k_{s+1}+1}{j_{s+1}+1}-\frac{k_s+1}{j_s+1}\right)B_s \\
\pagebreak[3]
& = & \left(\frac{k_s+1}{j_s-i_s}-\frac{k_{s+1}+1}{j_{s+1}+1}\right)(A_s-B_s)
\quad > \quad 0.
\end{eqnarray*}
Hence $F_0> F_r>0$, which contradicts (\ref{eq:F_0<0}).\QED

\begin{remark} \rm We have not worked out in general which values of $\nu(I_p)<\nu_0$ are forbidden, but point out that the symmetry of $H(R/I)$ restricts the possibilities of small values of $\nu(I_p)$ just as it rules out large values. For example, if $g=4$, then $\nu_0=(p+1)^2$. Let us see what happens in case $p=4$ and $\nu(I_4)<25$. There is no apparent Hilbert
function obstruction to values $\nu(I_4)\le15$, but in the range $15<
\nu(I_4) <25$ one will have $10<h_j=h_4<20$ growing to $h_{j+1}=20$,
and this suggests that arguments along the lines we have been giving
will still be effective. Since $a_j\ge j+2$ (notation of
 equation~(\ref{eq:h_j})) is only possible if $j=4$,
and $h_4={6\choose4}=15$ can grow to $h_5=20<{\binom75}$, such arguments
can only succeed for $11\le h_j\le 14$. In this case we are looking
for a gain $h_{j+1}-h_j\ge6$, so the expansion of $h_j$ will have to have at
least six terms of form $\binom{i+1}i$; this is clearly impossible for
$h_j\le 26$. We conclude that the values $21\le\nu(I_4)\le24$ do not
occur. If $16\le\nu(I_4)\le20$, then $H(R/I)$ fails to be unimodal.
There are no known examples of such sequences for $g=4$ at all, and
some evidence that they may not be possible. One can interpret our
arguments as ruling out ``extreme'' failure of unimodality, leaving a
grey zone of ``mild'' failure of unimodality open for more
investigation.\end{remark}

\begin{remark} It has been suggested to us that there is the possibility of a very short
and elegant argument, at least for the main result that $\nu(I_p)\le \nu_0$, based on 
the behavior of the combinatorial functions $f(x)=x^{\langle n \rangle}$ and $g(x)=
h^{\langle x \rangle }$. The key to this argument rests on the observation (for which we
do not know a proof) that $g$ is non-increasing. \end{remark}

\section{Resolutions and additional conjectures}

\hskip1\parindent We now give an example to show that Conjecture \ref{conj:1} can not be extended to bound the number of generators in all degrees. 

\begin{example} \rm Let 
$$
I=(x_1^2,\  x_1x_2x_3+x_3^2x_4,\  x_3^3,\  x_1x_3^2,\  x_2^4,\  x_4^4,\  x_1x_4^3,\  x_1x_2x_4^2+x_3x_4^3,\  x_2^3x_3^2,\  x_2^3x_4^3).
$$
 This is just the ideal quotient $(x_1^2, x_2^4, x_3^3, x_4^4)\link(x_1x_2-x_3x_4)$. Then $R/I$ is a Gorenstein artin algebra with $h$-vector $(1,4,9,13,13,9,4,1)$ and betti sequence $(1,10, 18, 10,1)$, whereas the $h$-vector for an extremal Gorenstein algebra of codimension four and initial degree two is $(1,4,1)$ and the betti sequence is $(1,9,16,9,1)$ (notice in particular that the multiplicity $e(R/I)=54$ is far greater than the minimal value of six). The graded structure of the minimal free resolution is given by the MACAULAY diagram in which the $(i,j)$ entry (starting with $(0,0)$ in the top left left hand corner) represents the rank of $R(-i-j)$ in the $j$th term of the resolution.
$$\begin{array}{ccccc}
      1   &  -  &   -  &   -  &   - \\
      -   &  1  &   -  &   -   &  - \\
      -   &  3  &   4   &  1   &  - \\
      -   &  4  &   5  &    1 &  - \\
      -   &  1   &  5  &   4   &  - \\
      -   &  1   &  4  &   3   &  - \\
      -   &  -   &  -  &   1 &   - \\
      -   &  -   &  -  &  -   &  1 
\end{array}
$$ 
\end{example}

\begin{remark} \rm If $I$ is generated only in degree $p$, and $g=4$, and $R/J$ is an extremal algebra of the same codimension and initial degree, then the theorem shows that $\b_1(R/I)\le\b_1(R/J)$; symmetry of the minimal free resolution and $\sum_{i=0}^4 (-1)^i\b_i(R/I)=0$ then imply that $\b_i(R/I)\le \b_i(R/J)$ for all $i$.
\end{remark}

It seems likely that in any codimension if $R/I$ has a pure minimal free resolution, then $\b_i(R/I)\le\b_i(R/J)$ for all $i$. This stronger estimate is true in codimension three, without any purity assumption, thanks to the Buchsbaum-Eisenbud structure theorem. For initial degree $p$, the extremal Gorenstein algebra of codimension three has $\b_1=\b_2=2p+1$.

\begin{theorem}\label{codim3} Let $R$ be a polynomial ring over a field and $I$ be a homogeneous Gorenstein ideal of height three. If $p$ is the initial degree of $I$, then $\nu(I)\le 2p+1$ and $\b_2(R/I)\le 2p+1$. 
\end{theorem}
\demo We may assume, without loss of generality, that $R=k[x_1,x_2,x_3]$ and $A=R/I$ is artinian and local with $\socle(A)=A_\sigma$. Then by \cite{BE} the minimal free resolution of $A$ has the form 
$$
0\to R(-\sigma-3)\to \oplus_{j=1}^\nu R(-n_j) \stackrel{Y}{\longrightarrow} \oplus_{i=1}^\nu R(-m_i)\to R,
$$ where $Y$ is an alternating matrix, and the generators $f_1,\dots,f_\nu $ of $I$ are the maximal pfaffians of $Y$. We may take $p=m_1\le m_2\le\dots\le m_\nu$ and $n_1\le n_2\le\dots\le n_\nu=\sigma+3-p$. Let $d_{ij}=\deg(y_{ij})=
m_i-n_j$; notice that $y_{ij}=0$ by minimality of the resolution if $d_{ij}\le0$. Each generator is a sum of monomials of degree $(\nu-1)/2$ in the $y_{ij}$ (evident from the Laplace expansion for pfaffians); and any such term that is non-zero has degree at least $(\nu-1)/2$ in $R$. Hence if the theorem fails every generator has degree at least $(2p+1)/2>p$, which is a contradiction since at least one minimal generator has degree $p$.  \QED

\section{Acknowledgements}

\hskip1\parindent
Investigation of this problem first begun when the first author was on sabbatical at Rutgers University, and continued during a visit to the Instituto Polit\'ecnico Nacional; he would like to thank both institutions for their support and hospitality.

\end{document}